\documentclass[reqno]{amsart}

\usepackage{lmodern}
\usepackage[T1]{fontenc}
\usepackage[utf8]{inputenc}
\usepackage[english]{babel}
\usepackage{microtype} 

\usepackage{amsmath,amssymb,amsthm,mathrsfs,latexsym,mathtools,mathdots,enumerate,tikz,bm,url}
\usepackage[enableskew,vcentermath]{youngtab}
\usepackage[centertableaux]{ytableau}

\usetikzlibrary{arrows,positioning}
\usepackage[capitalize]{cleveref}

\usepackage{todonotes}
\presetkeys%
    {todonotes}%
    {inline,backgroundcolor=yellow}{}

\usepackage{graphicx}

\usepackage{ifpdf}
\graphicspath{{./figures/}}
\usepackage{epstopdf}
\DeclareGraphicsExtensions{.png,.jpg,.eps,.epsf}

\newtheorem{theorem}{Theorem}
\newtheorem{proposition}[theorem]{Proposition}

\theoremstyle{definition}

\newtheorem{example}[theorem]{Example}
\newtheorem{remark}[theorem]{Remark}

\newcommand{\defin}[1]{\emph{#1}}

\newcommand{\setZ}{\mathbb{Z}}

\newcommand{\setR}{\mathbb{R}}

\newcommand{\xvec}{\mathbf{x}}

\newcommand{\schurS}{\mathrm{s}}

\newcommand{\oPoly}{\mathcal{O}}
\newcommand{\polyG}{\mathcal{G}}
\newcommand{\polyP}{\mathcal{P}}

\newcommand{\ehr}{\mathrm{ehr}}

\DeclareMathOperator{\hook}{hook}

\title{Polytopes and large counterexamples}
\author{Per Alexandersson}
\address{Dept. of Mathematics, University of Pennsylvania. Philadelphia, PA}
\email{per.w.alexandersson@gmail.com}
\subjclass[2010]{52B12, 05E10}

\begin{document}

\begin{abstract}
In this short note, we give large counterexamples to natural questions about certain order polytopes,
in particular, Gelfand--Tsetlin polytopes. 
Several of the counterexamples are too large to be discovered via a brute-force computer search.

We also show that the multiset of hooks in a Young diagram is not enough information 
to determine the Ehrhart polynomial for an associated order polytope.
This is somewhat counter-intuitive to the fact that the 
multiset of hooks always determine the leading coefficient of the Ehrhart polynomial.
\end{abstract}

\maketitle

\section{Introduction}

It is nowadays easy to make conjectures based on extensive computer evidence. 
This note gives several examples where the brute-force approach yield strong support for natural conjectures,
but turns out to fail in higher dimensions.

There is a long-standing open conjecture regarding Kostka numbers, with plenty of computational support.
It is conjectured, \cite{King04stretched} that the map $k \to K_{k \lambda / k\mu, k w}$ 
is a polynomial with non-negative coefficients.
In particular, this conjecture implies that the Ehrhart polynomial of the Birkhoff polytopes 
(defined as the convex hull of all $n\times n$ permutation matrices) has non-negative coefficients ---
also an open problem.

The Kostka coefficients can be computed by counting the lattice points in 
the intersection of a hyperplane defined by $w$ and the Gelfand--Tsetlin polytope $\polyG_{\lambda/\mu}$.
It is therefore natural to study the Ehrhart polynomial of other restrictions of $\polyG_{\lambda/\mu}$.
In \cref{ex:gtface}, we give an example of a face of a certain $\polyG_{\lambda/\mu}$
with a negative coefficient. The example is too large to be discovered by brute force computer search,
suggesting that eventual counterexamples to the Koskta numbers conjecture is of a similar nature.

\medskip

One can show that the Gelfand--Tsetlin polytopes $\polyG_{\lambda/\mu}$
have the integer decomposition property (IDP), and  that certain intersections with 
hyperplanes also give rise to this property, see \cite{Alexandersson2016GTPoly}.
We present a large (carefully constructed) counterexample to a natural conjecture in this setting.

A few recent papers present conjectures regarding the IDP and Ehrhart coefficients,
see \cite{Braun2016} and \cite{Hibi2016}.
The latter regards \emph{lecture hall polytopes}, 
which are closely related to the polytopes we consider below in \cref{ex:partitionpolytope},
suggesting that possible counterexamples to conjectures in \cite{Hibi2016}
only occur in high dimensions.

\medskip

Specifying the multiset of hooks gives very strong restrictions on
a Young diagram, and therefore, a related order polytope.
We show that the multiset of hooks of a Young diagram \emph{does not} determine the Ehrhart polynomial
of the order polytope --- even though the multiset determines the leading coefficient.
A similar example is provided  for trees.

\subsection{Order polytopes}

We follow the terminology in \cite{Stanley86TwoPosetPolytopes,StanleyEC1}. 
Given a partial order $P$ on the set $[n]$, we associate a polytope $\oPoly(P)$, called the \defin{order polytope} of $P$.
It is the polytope in $\setR^n$ defined via the inequalities
\begin{align}
x_i \leq x_j \text{ if } i <_P j \qquad \text{and} \qquad 0 \leq x_i \leq 1 \text{ for } 1\leq i \leq n.
\end{align}

There is a slight generalization of order polytopes, called \defin{marked order polytopes}, see \emph{e.g.,} \cite{ArdilaBS11},
where we replace the conditions $0 \leq x_i \leq 1$ with $a_i \leq x_i \leq b_i$ for $1\leq i \leq n$,
where the $a_i$ and $b_i$ are fixed constants.
We will mostly be concerned with a special case of 
marked order polytopes, the Gelfand--Tsetlin polytopes, which we define below.

\medskip

One property worth mentioning is that all marked order polytopes with \emph{integer boundary conditions} 
$a_i$ and $b_i$, are integral, \emph{i.e.}, all vertices are integer lattice points. 
In fact, integral marked order polytopes admit a \defin{unimodular triangulation} --- a triangulation 
into simplices with integer vertices, all having normalized volume $1$.
To sketch a proof of the latter statement, suppose that the polytope is in $\setR^n$ and consider the hyperplanes 
\begin{align}
x_i=m \text{ and } x_i = x_{j} + m \text{ for all } 1 \leq i, j \leq n \text{ and } m \in \setZ. \label{eq:hyperplanes}
\end{align}
The cube with vertices in $\{0,1\}^n$ is partitioned into $n!$ unimodular simplices by these hyperplanes,
implying that \eqref{eq:hyperplanes} determines a unimodular triangulation of  $\setR^n$.
Finally note that the defining hyperplanes of a marked 
order polytope is a subset of the hyperplanes in \eqref{eq:hyperplanes}.

\subsection{Gelfand--Tsetlin polytopes}

Gelfand--Tsetlin polytopes were first introduced in \cite{GelfandTsetlin50}.
Consider a parallelogram arrangement of non-negative numbers,
\[
\begin{array}{cccccccccccccc}
x^m_{1} & & x^m_{2} & & \cdots & & \cdots & & x^m_{n} \\
 & \ddots & & \ddots &  & &   & & & \ddots  \\
  &  &   x^1_{1} &  & x^1_{2} & & \cdots & & \cdots & & x^1_{n} \\
   &  &    & x^0_{1} &    & x^0_{2} & & \cdots & & \cdots & & x^0_{n}
\end{array}
\]
where the entries must satisfy the inequalities
\begin{equation}
x^{i+1}_{j} \geq x^i_{j} \text{ and } x^i_{j} \geq x^{i+1}_{j+1} \label{eq:gtinequalities}
\end{equation}
for all values of $i$, $j$ where the indexing is defined.
Note that horizontal rows and down-right diagonals are weakly decreasing,
while down-left diagonals are weakly increasing. 

Let $\lambda \supset \mu$ be partitions and fix the top and bottom rows in the 
diagram to be $\lambda$ and $\mu$, \emph{i.e.,} $x^m_i = \lambda_i$ and $x^0_i = \mu_i$ for $1\leq i \leq n$.
These conditions together with the inequalities in \cref{eq:gtinequalities} defines a polytope $\polyG_{\lambda/\mu}$.
Note that this polytope also depends on the choice of $m$.
One can show that lattice points in $\polyG_{\lambda/\mu}$ are in bijection with semi-standard Young tableaux
with shew shape $\lambda/\mu$ and entries in $[m]$, see \emph{e.g.}, \cite{StanleyEC2}.

\medskip 
In the case $\mu$ is the empty partition, many of the entries in the parallel arrangement are forced to be $0$.
It suffices to keep track of the triangular pattern
\[
\begin{array}{cccccccccccccc}
x^m_{1} & & x^m_{2} & & & \cdots & & x^m_{n} \\
 & \ddots & & \ddots &    & \iddots \\
  &  &   x^2_{1} &  & x^2_{2}  \\
   &  &    & x^1_{1} 
\end{array}
\]
where the top row is fixed to be $\lambda$. We denote this polytope $\polyG_\lambda$.

\subsection{Ehrhart polynomial}

Given a polytope $\polyP \subset \setR^n$ with integer vertices, the \emph{Ehrhart polynomial} $\ehr_\polyP(n)$
is the function defined as
\begin{align*}
 \ehr_\polyP(n) = \#\left(n \polyP \cap \setZ^n \right).
\end{align*}
In other words, $\ehr_\polyP(n)$ counts the number of lattice points in the $n$th dilate of $\polyP$.
Ehrhart showed that this is a polynomial in $n$.

\medskip 

We will use the notation $\ehr_{\lambda/\mu}(n)$ or simply $\ehr_{\lambda}(n)$ to denote the Ehrhart polynomials 
associated with $\polyG_{\lambda/\mu}$ and $\polyG_\lambda$, respectively.
One can show that 
\begin{align*}
\ehr_{\lambda/\mu}(n) = \schurS_{n\lambda/n\mu}(1^{m})
\end{align*}
where $\schurS_{\lambda/\mu}(\xvec)$ is the skew Schur function and $1^m = (1,1,\dotsc,1)$ with $m$ ones
and $n\lambda$ is interpreted as elementwise multiplication.
Remember, $m+1$ is the number of rows in the parallelogram arrangement defining the Gelfand--Tsetlin polytope.

\section{Order polytopes and faces of Gelfand--Tsetlin polytopes}

Now we are ready to give the first counterexample to a reasonable conjecture.
An exhaustive search on the computer verifies that all order polytopes 
with dimension $\leq 7$ have no negative coefficients in their Ehrhart polynomials.

\begin{example}[Negative Ehrhart coefficients]\label{ex:orderpoly}
There are order polytopes with negative Ehrhart coefficients.
Consider the following polytope, given by the inequalities $0 \leq z \leq x_i \leq 1$ for $i=1,\dotsc,\ell$:
\begin{figure}[!ht]
\begin{tikzpicture}[scale=0.8]
\node (top) at (0,0) {$1$};
\node (v1) at (-4,-1) {$x_1$};
\node (v2) at (-2,-1) {$x_2$};
\node (vd) at (0,-1) {$\dots$};
\node (v3) at (2,-1) {$x_{\ell-1}$};
\node (v4) at (4,-1) {$x_\ell$};
\node (bot1) at (0,-2) {$z$};
\node (bot0) at (0,-3) {$0$};
\draw  (top) -- (v1)--(bot1);
\draw  (top) -- (v2)--(bot1);
\draw  (top) -- (v3)--(bot1);
\draw  (top) -- (v4)--(bot1);
\draw  (bot1) -- (bot0);
\end{tikzpicture}
\caption{For $\ell=20$, there are negative coefficients.}
\end{figure}
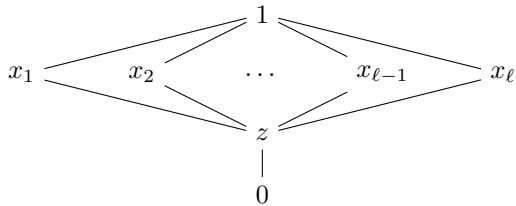
It is fairly straightforward to show that this order polytope has Ehrhart 
polynomial given by $\ehr(n) = \sum_{j=1}^{n+1} j^\ell$.
This expression can be evaluated using a computer algebra system and for $\ell=20$ we get
\[
 \ehr(n) = \frac{1}{6930}\left(6930-3528231n+1316700 n^2 + 32027050 n^3 + \dotsb \right).
\]
For $\ell < 20$, all coefficients are non-negative.
Note that even using state-of-the-art lattice point counting software such as \texttt{lattE}
to determine the Ehrhart polynomial of this order polytope is out of reach due to the high dimension.

\noindent
\emph{Note:} This construction more or less appear in \cite[Figure 3.87]{StanleyEC1}
and should be attributed to R.\ Stanley \cite{StanleyMOEhrhart}.
\end{example}

It is easy to see from the Weyl character formula that for a non-skew shape $\lambda$,
the Ehrhart polynomial $\ehr_{\lambda}(n)$ has non-negative coefficients.
A natural question is then to consider \emph{faces} of Gelfand--Tsetlin polytopes.

Faces of Gelfand--Tsetlin polytopes are obtained by 
forcing some of the inequalities in \cref{eq:gtinequalities} to be equalities.
Such faces appear in many places in representation theory and algebraic combinatorics.
For example, certain faces are responsible for generating key polynomials, see \cite{Kiritchenko2010},
and they are also in correspondence with so called $rc$-graphs and pipe dreams, see \cite{Kogan2005}.
Flagged Schur polynomials can be computed as a sum over lattice points in certain faces of Gelfand--Tsetlin polytopes.

We can use the previous example to construct a face of a Gelfand--Tsetlin polytope $\polyG_\lambda$,
such that a certain face has a negative coefficient in the Ehrhart polynomial.

\begin{example}\label{ex:gtface}
We can construct a one-parameter family of faces of Gelfand--Tsetlin polytopes with top row given by $\lambda = (1^\ell, 0^{\ell+1})$:
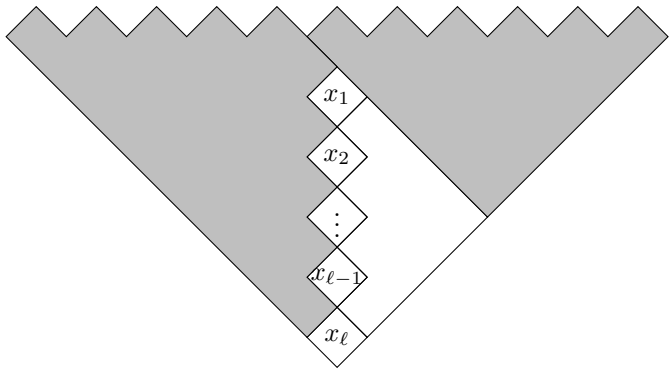
\[
\begin{tikzpicture}[scale=0.4]
\draw[black] (12,-2)--(13,-3)--(14,-2)--(13,-1)--(12,-2)--cycle;
\draw[black] (13,-9)--(14,-10)--(15,-9)--(16,-8)--(17,-7)--(18,-6)--(17,-5)--(16,-4)--(15,-3)--(14,-2)--(13,-3)--(14,-4)--(13,-5)--(14,-6)--(13,-7)--(14,-8)--(13,-9)--cycle;
\draw[black] (12,-4)--(13,-5)--(14,-4)--(13,-3)--(12,-4)--cycle;
\draw[black] (12,-6)--(13,-7)--(14,-6)--(13,-5)--(12,-6)--cycle;
\draw[black] (12,-8)--(13,-9)--(14,-8)--(13,-7)--(12,-8)--cycle;
\draw[black] (12,-10)--(13,-11)--(14,-10)--(13,-9)--(12,-10)--cycle;
\draw[color=black,fill=lightgray] (18,-6)--(19,-5)--(20,-4)--(21,-3)--(22,-2)--(23,-1)--(24,0)--(23,1)--(22,0)--(21,1)--(20,0)--(19,1)--(18,0)--(17,1)--(16,0)--(15,1)--(14,0)--(13,1)--(12,0)--(13,-1)--(14,-2)--(15,-3)--(16,-4)--(17,-5)--(18,-6)--cycle;
\draw[color=black,fill=lightgray] (12,-10)--(13,-9)--(12,-8)--(13,-7)--(12,-6)--(13,-5)--(12,-4)--(13,-3)--(12,-2)--(13,-1)--(12,0)--(11,1)--(10,0)--(9,1)--(8,0)--(7,1)--(6,0)--(5,1)--(4,0)--(3,1)--(2,0)--(3,-1)--(4,-2)--(5,-3)--(6,-4)--(7,-5)--(8,-6)--(9,-7)--(10,-8)--(11,-9)--(12,-10)--cycle;
\node at (3,0) {$ $};
\node at (5,0) {$ $};
\node at (7,0) {$ $};
\node at (9,0) {$ $};
\node at (11,0) {$ $};
\node at (13,0) {$ $};
\node at (15,0) {$ $};
\node at (17,0) {$ $};
\node at (19,0) {$ $};
\node at (21,0) {$ $};
\node at (23,0) {$ $};
\node at (4,-1) {$ $};
\node at (6,-1) {$ $};
\node at (8,-1) {$ $};
\node at (10,-1) {$ $};
\node at (12,-1) {$ $};
\node at (14,-1) {$ $};
\node at (16,-1) {$ $};
\node at (18,-1) {$ $};
\node at (20,-1) {$ $};
\node at (22,-1) {$ $};
\node at (5,-2) {$ $};
\node at (7,-2) {$ $};
\node at (9,-2) {$ $};
\node at (11,-2) {$ $};
\node at (13,-2) {$x_1$};
\node at (15,-2) {$ $};
\node at (17,-2) {$ $};
\node at (19,-2) {$ $};
\node at (21,-2) {$ $};
\node at (6,-3) {$ $};
\node at (8,-3) {$ $};
\node at (10,-3) {$ $};
\node at (12,-3) {$ $};
\node at (14,-3) {$ $};
\node at (16,-3) {$ $};
\node at (18,-3) {$ $};
\node at (20,-3) {$ $};
\node at (7,-4) {$ $};
\node at (9,-4) {$ $};
\node at (11,-4) {$ $};
\node at (13,-4) {$x_2$};
\node at (15,-4) {$ $};
\node at (17,-4) {$ $};
\node at (19,-4) {$ $};
\node at (8,-5) {$ $};
\node at (10,-5) {$ $};
\node at (12,-5) {$ $};
\node at (14,-5) {$ $};
\node at (16,-5) {$ $};
\node at (18,-5) {$ $};
\node at (9,-6) {$ $};
\node at (11,-6) {$ $};
\node at (13,-6) {$\mathscr{\vdots}$};
\node at (15,-6) {$ $};
\node at (17,-6) {$ $};
\node at (10,-7) {$ $};
\node at (12,-7) {$ $};
\node at (14,-7) {$ $};
\node at (16,-7) {$ $};
\node at (11,-8) {$ $};
\node at (13,-8) {$x_{\ell-1}$};
\node at (15,-8) {$ $};
\node at (12,-9) {$ $};
\node at (14,-9) {$ $};
\node at (13,-10) {$x_\ell$};
\end{tikzpicture}
\]
There are two shaded regions and one larger unshaded region, where all entries within each respective region are all forced to be equal.
These conditions determine a face. Note that all entries within the leftmost shaded region are all equal to $1$,
while all variables in the rightmost region are equal to $0$.
Lattice points in this face are in clear bijection with lattice points in the order polytope in \cref{ex:orderpoly},
so for $\ell = 20$, the Ehrhart polynomial has a negative coefficient.
\end{example}

\section{Partition polytope and a large GT counterexample}

As in the previous section, we first construct a counterexample in a simpler setting,
which is then transferred to a statement about Gelfand--Tsetlin polytopes.
The property we consider in this section is the integer decomposition property.

An integral polytope $\polyP \subset \setR^d$ is said to have the \defin{integer decomposition property} (IDP) 
if for every positive integer $k$ and $\xvec \in k \polyP \cap \setZ^d$,
we can find $\xvec_1,\xvec_2,\dots,\xvec_k \in \polyP \cap \setZ^d$ such that $\xvec_1 + \dots + \xvec_k = \xvec$.
One can show the following series of implications for a polytope $\polyP$;
\[
 \text{$\polyP$ has a unimodular triangulation} \Longrightarrow \text{$\polyP$ has the IDP} \Longrightarrow \text{$\polyP$ is integral}.
\]

Let $\polyP_{a,b}$ denote the convex cull of all partitions of $a$ with $b$ parts, seen as integer points in $\setR^b$.
Computer experiments suggests that $\polyP_{a,b}$ has the IDP if $a \leq 18$ and  $b \leq 9$.
However, we have the following counterexample, which was pointed out by Robert Davis to the author via private communication:

\begin{example}\label{ex:partitionpolytope}
The polytope $\polyP_{18,9}$ does not have the IDP.
Consider the lattice point $p = (6,6,6,6,4,4,2,1,1)$. It is straightforward to verify that
\begin{align}\label{eq:4partitions}
\begin{split}
p =& \frac12(4,4,4,4,1,1,0,0,0) +\frac12 (3,3,3,3,3,3,0,0,0) + \\
&\frac12 (3,3,3,3,2,2,2,0,0) + \frac12 (2,2,2,2,2,2,2,2,2)
\end{split}
\end{align}
so it is clear that $p \in 2\cdot \polyP_{18,9}$. However, an exhaustive search on the computer
shows that $p$ cannot be expressed as a sum of two lattice points in $\polyP_{18,9}$.
\end{example}
The point $p$ above is not the only point in $2\cdot \polyP_{18,9}$ that cannot be expressed 
as a sum of two integer points in $\polyP_{18,9}$ --- another such point is $(6,6,5,5,5,4,2,2,1)$.

Note that each row in a Gelfand--Tsetlin pattern is a partition. 
Let us add an extra restriction on Gelfand--Tsetlin polytopes, and let $\polyG_{\lambda/\mu,w}$
denote the polytope consisting of all points satisfying \cref{eq:gtinequalities},
with the additional condition that the entries in row $i$ sum up to $w_i$.
Here, $w$ is an integer vector with $m$ rows.
The polytopes $\polyG_{\lambda/\mu,w}$ are in general no longer integral, 
a fact first pointed out in \cite{King04stretched}, and further extended in \cite{Loera04,LoeraM06,Alexandersson2016GTPoly}.
\medskip

Since Gelfand--Tsetlin polytopes have nice properties with respect to unimodular triangulations and the IDP,
it is natural to ask if similar properties hold for $\polyG_{\lambda/\mu,w}$.
In fact, in \cite{Alexandersson2016GTPoly} we show that $\polyG_{\lambda/\mu,1^m}$
always has a unimodular triangulation whenever it is integral.

Can we do something in the non-integral case?
By using \cref{ex:partitionpolytope} as a starting point, we can show the following statement:
\begin{proposition}
Let $\polyG^\ast_{\lambda/\mu,w}$ be the polytope defined as the 
convex hull of the integer points in $\polyG^\ast_{\lambda/\mu,w}$.
Then there is a $\polyG^\ast_{\lambda/\mu,w}$ which does not satisfy the integer decomposition property.
\end{proposition}
We have not been able to find a $\polyG^\ast_{\lambda/\mu,w}$ which does not have the IDP 
by brute force computer search.
The following is the smallest known counterexample which has been carefully constructed:
\begin{example}
Let
\begin{align*}
 \lambda &= (4,4,4,4,3,3,2,2,2), \quad  \mu = (2,0,0,0,0,0,0,0,0) \\
 w &= (2, 4, 5, 6, 7, 8, 9, 10, 11, 12, 14, 16, 18, 20, 22, 23, 25, 27, 28). 
\end{align*}
Then $\polyG^\ast_{\lambda/\mu,w}$ does not have the IDP. Consider the following GT-pattern $G$:
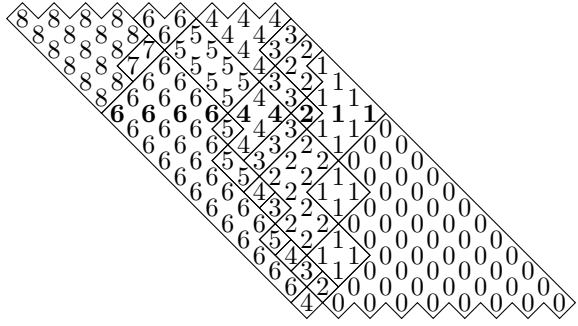
\begin{align*}
\begin{tikzpicture}[scale=0.21]
\draw[black] (20,-18)--(21,-19)--(22,-18)--(21,-17)--(20,-18)--cycle;
\draw[black] (14,0)--(13,1)--(12,0)--(11,1)--(10,0)--(11,-1)--(12,-2)--(13,-1)--(14,0)--cycle;
\draw[black] (18,0)--(17,1)--(16,0)--(15,1)--(14,0)--(15,-1)--(16,-2)--(17,-3)--(18,-4)--(19,-3)--(18,-2)--(19,-1)--(20,0)--(19,1)--(18,0)--cycle;
\draw[black] (9,-5)--(10,-4)--(9,-3)--(10,-2)--(11,-1)--(10,0)--(9,1)--(8,0)--(7,1)--(6,0)--(5,1)--(4,0)--(3,1)--(2,0)--(3,-1)--(4,-2)--(5,-3)--(6,-4)--(7,-5)--(8,-6)--(9,-5)--cycle;
\draw[black] (22,-18)--(23,-19)--(24,-18)--(25,-19)--(26,-18)--(27,-19)--(28,-18)--(29,-19)--(30,-18)--(31,-19)--(32,-18)--(33,-19)--(34,-18)--(35,-19)--(36,-18)--(37,-19)--(38,-18)--(37,-17)--(36,-16)--(35,-15)--(34,-14)--(33,-13)--(32,-12)--(31,-11)--(30,-10)--(29,-9)--(28,-8)--(27,-7)--(26,-6)--(25,-7)--(24,-8)--(23,-9)--(24,-10)--(25,-11)--(24,-12)--(23,-13)--(24,-14)--(25,-15)--(24,-16)--(23,-17)--(22,-18)--cycle;
\draw[black] (17,-5)--(18,-4)--(17,-3)--(16,-2)--(15,-1)--(14,0)--(13,-1)--(12,-2)--(13,-3)--(14,-4)--(15,-5)--(16,-6)--(17,-5)--cycle;
\draw[black] (18,-2)--(19,-3)--(20,-2)--(21,-1)--(20,0)--(19,-1)--(18,-2)--cycle;
\draw[black] (9,-3)--(10,-4)--(11,-3)--(12,-2)--(11,-1)--(10,-2)--(9,-3)--cycle;
\draw[black] (21,-3)--(22,-2)--(21,-1)--(20,-2)--(19,-3)--(20,-4)--(21,-5)--(22,-4)--(21,-3)--cycle;
\draw[black] (18,-16)--(19,-17)--(20,-18)--(21,-17)--(20,-16)--(19,-15)--(18,-14)--(19,-13)--(18,-12)--(17,-11)--(16,-10)--(15,-9)--(16,-8)--(15,-7)--(16,-6)--(15,-5)--(14,-4)--(13,-3)--(12,-2)--(11,-3)--(10,-4)--(9,-5)--(8,-6)--(9,-7)--(10,-8)--(11,-9)--(12,-10)--(13,-11)--(14,-12)--(15,-13)--(16,-14)--(17,-15)--(18,-16)--cycle;
\draw[black] (24,-8)--(25,-7)--(26,-6)--(25,-5)--(24,-4)--(23,-3)--(22,-2)--(21,-3)--(22,-4)--(21,-5)--(22,-6)--(21,-7)--(22,-8)--(23,-9)--(24,-8)--cycle;
\draw[black] (20,-4)--(19,-3)--(18,-4)--(19,-5)--(20,-6)--(21,-5)--(20,-4)--cycle;
\draw[black] (18,-8)--(19,-7)--(20,-6)--(19,-5)--(18,-4)--(17,-5)--(16,-6)--(17,-7)--(16,-8)--(17,-9)--(18,-8)--cycle;
\draw[black] (20,-6)--(21,-7)--(22,-6)--(21,-5)--(20,-6)--cycle;
\draw[black] (15,-7)--(16,-8)--(17,-7)--(16,-6)--(15,-7)--cycle;
\draw[black] (19,-9)--(20,-8)--(21,-7)--(20,-6)--(19,-7)--(18,-8)--(17,-9)--(18,-10)--(19,-9)--cycle;
\draw[black] (19,-13)--(20,-14)--(21,-15)--(22,-14)--(23,-13)--(22,-12)--(21,-11)--(22,-10)--(23,-9)--(22,-8)--(21,-7)--(20,-8)--(19,-9)--(18,-10)--(19,-11)--(20,-12)--(19,-13)--cycle;
\draw[black] (17,-9)--(16,-8)--(15,-9)--(16,-10)--(17,-11)--(18,-10)--(17,-9)--cycle;
\draw[black] (24,-12)--(25,-11)--(24,-10)--(23,-9)--(22,-10)--(21,-11)--(22,-12)--(23,-13)--(24,-12)--cycle;
\draw[black] (17,-11)--(18,-12)--(19,-11)--(18,-10)--(17,-11)--cycle;
\draw[black] (18,-12)--(19,-13)--(20,-12)--(19,-11)--(18,-12)--cycle;
\draw[black] (18,-14)--(19,-15)--(20,-14)--(19,-13)--(18,-14)--cycle;
\draw[black] (24,-16)--(25,-15)--(24,-14)--(23,-13)--(22,-14)--(21,-15)--(22,-16)--(23,-17)--(24,-16)--cycle;
\draw[black] (19,-15)--(20,-16)--(21,-15)--(20,-14)--(19,-15)--cycle;
\draw[black] (20,-16)--(21,-17)--(22,-16)--(21,-15)--(20,-16)--cycle;
\draw[black] (21,-17)--(22,-18)--(23,-17)--(22,-16)--(21,-17)--cycle;
\node at (3,0) {$8$};
\node at (5,0) {$8$};
\node at (7,0) {$8$};
\node at (9,0) {$8$};
\node at (11,0) {$6$};
\node at (13,0) {$6$};
\node at (15,0) {$4$};
\node at (17,0) {$4$};
\node at (19,0) {$4$};
\node at (4,-1) {$8$};
\node at (6,-1) {$8$};
\node at (8,-1) {$8$};
\node at (10,-1) {$8$};
\node at (12,-1) {$6$};
\node at (14,-1) {$5$};
\node at (16,-1) {$4$};
\node at (18,-1) {$4$};
\node at (20,-1) {$3$};
\node at (5,-2) {$8$};
\node at (7,-2) {$8$};
\node at (9,-2) {$8$};
\node at (11,-2) {$7$};
\node at (13,-2) {$5$};
\node at (15,-2) {$5$};
\node at (17,-2) {$4$};
\node at (19,-2) {$3$};
\node at (21,-2) {$2$};
\node at (6,-3) {$8$};
\node at (8,-3) {$8$};
\node at (10,-3) {$7$};
\node at (12,-3) {$6$};
\node at (14,-3) {$5$};
\node at (16,-3) {$5$};
\node at (18,-3) {$4$};
\node at (20,-3) {$2$};
\node at (22,-3) {$1$};
\node at (7,-4) {$8$};
\node at (9,-4) {$8$};
\node at (11,-4) {$6$};
\node at (13,-4) {$6$};
\node at (15,-4) {$5$};
\node at (17,-4) {$5$};
\node at (19,-4) {$3$};
\node at (21,-4) {$2$};
\node at (23,-4) {$1$};
\node at (8,-5) {$8$};
\node at (10,-5) {$6$};
\node at (12,-5) {$6$};
\node at (14,-5) {$6$};
\node at (16,-5) {$5$};
\node at (18,-5) {$4$};
\node at (20,-5) {$3$};
\node at (22,-5) {$1$};
\node at (24,-5) {$1$};
\node at (9,-6) {$\mathbf 6$};
\node at (11,-6) {$\mathbf 6$};
\node at (13,-6) {$\mathbf 6$};
\node at (15,-6) {$\mathbf 6$};
\node at (17,-6) {$\mathbf 4$};
\node at (19,-6) {$\mathbf 4$};
\node at (21,-6) {$\mathbf 2$};
\node at (23,-6) {$\mathbf 1$};
\node at (25,-6) {$\mathbf 1$};
\node at (10,-7) {$6$};
\node at (12,-7) {$6$};
\node at (14,-7) {$6$};
\node at (16,-7) {$5$};
\node at (18,-7) {$4$};
\node at (20,-7) {$3$};
\node at (22,-7) {$1$};
\node at (24,-7) {$1$};
\node at (26,-7) {$0$};
\node at (11,-8) {$6$};
\node at (13,-8) {$6$};
\node at (15,-8) {$6$};
\node at (17,-8) {$4$};
\node at (19,-8) {$3$};
\node at (21,-8) {$2$};
\node at (23,-8) {$1$};
\node at (25,-8) {$0$};
\node at (27,-8) {$0$};
\node at (12,-9) {$6$};
\node at (14,-9) {$6$};
\node at (16,-9) {$5$};
\node at (18,-9) {$3$};
\node at (20,-9) {$2$};
\node at (22,-9) {$2$};
\node at (24,-9) {$0$};
\node at (26,-9) {$0$};
\node at (28,-9) {$0$};
\node at (13,-10) {$6$};
\node at (15,-10) {$6$};
\node at (17,-10) {$5$};
\node at (19,-10) {$2$};
\node at (21,-10) {$2$};
\node at (23,-10) {$1$};
\node at (25,-10) {$0$};
\node at (27,-10) {$0$};
\node at (29,-10) {$0$};
\node at (14,-11) {$6$};
\node at (16,-11) {$6$};
\node at (18,-11) {$4$};
\node at (20,-11) {$2$};
\node at (22,-11) {$1$};
\node at (24,-11) {$1$};
\node at (26,-11) {$0$};
\node at (28,-11) {$0$};
\node at (30,-11) {$0$};
\node at (15,-12) {$6$};
\node at (17,-12) {$6$};
\node at (19,-12) {$3$};
\node at (21,-12) {$2$};
\node at (23,-12) {$1$};
\node at (25,-12) {$0$};
\node at (27,-12) {$0$};
\node at (29,-12) {$0$};
\node at (31,-12) {$0$};
\node at (16,-13) {$6$};
\node at (18,-13) {$6$};
\node at (20,-13) {$2$};
\node at (22,-13) {$2$};
\node at (24,-13) {$0$};
\node at (26,-13) {$0$};
\node at (28,-13) {$0$};
\node at (30,-13) {$0$};
\node at (32,-13) {$0$};
\node at (17,-14) {$6$};
\node at (19,-14) {$5$};
\node at (21,-14) {$2$};
\node at (23,-14) {$1$};
\node at (25,-14) {$0$};
\node at (27,-14) {$0$};
\node at (29,-14) {$0$};
\node at (31,-14) {$0$};
\node at (33,-14) {$0$};
\node at (18,-15) {$6$};
\node at (20,-15) {$4$};
\node at (22,-15) {$1$};
\node at (24,-15) {$1$};
\node at (26,-15) {$0$};
\node at (28,-15) {$0$};
\node at (30,-15) {$0$};
\node at (32,-15) {$0$};
\node at (34,-15) {$0$};
\node at (19,-16) {$6$};
\node at (21,-16) {$3$};
\node at (23,-16) {$1$};
\node at (25,-16) {$0$};
\node at (27,-16) {$0$};
\node at (29,-16) {$0$};
\node at (31,-16) {$0$};
\node at (33,-16) {$0$};
\node at (35,-16) {$0$};
\node at (20,-17) {$6$};
\node at (22,-17) {$2$};
\node at (24,-17) {$0$};
\node at (26,-17) {$0$};
\node at (28,-17) {$0$};
\node at (30,-17) {$0$};
\node at (32,-17) {$0$};
\node at (34,-17) {$0$};
\node at (36,-17) {$0$};
\node at (21,-18) {$4$};
\node at (23,-18) {$0$};
\node at (25,-18) {$0$};
\node at (27,-18) {$0$};
\node at (29,-18) {$0$};
\node at (31,-18) {$0$};
\node at (33,-18) {$0$};
\node at (35,-18) {$0$};
\node at (37,-18) {$0$};
\end{tikzpicture}
\end{align*}
We can verify that $G$ is lattice point in $2 \cdot \polyG^\ast_{\lambda/\mu,w}$,
since $G = \frac12(G_1 + G_2 + G_3 + G_4)$, where $G_1, \dotsc, G_4 \in \polyG^\ast_{\lambda/\mu,w}$ are given by
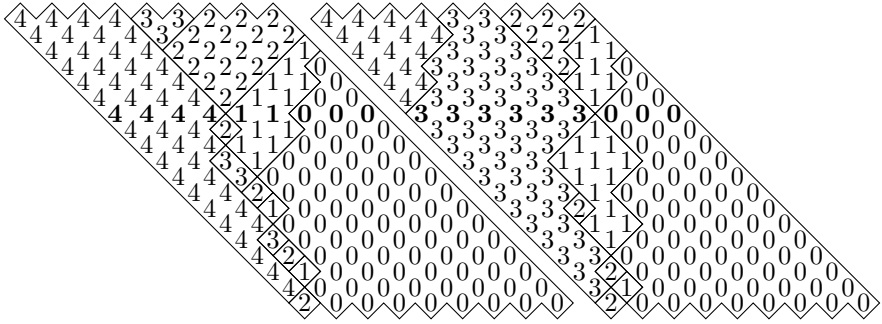
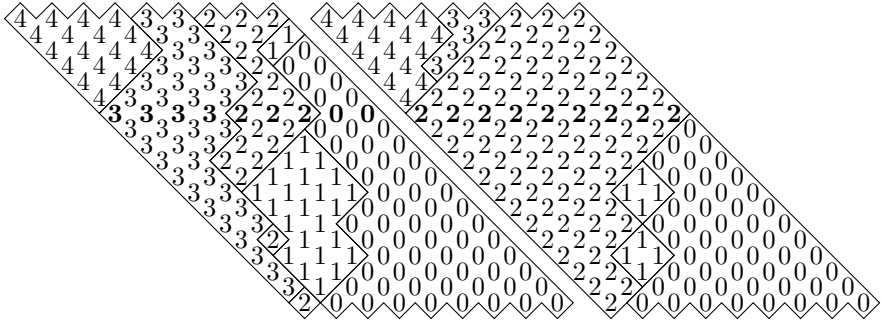
\[
\begin{tikzpicture}[scale=0.21]
\draw[black] (20,-18)--(21,-19)--(22,-18)--(21,-17)--(20,-18)--cycle;
\draw[black] (14,0)--(13,1)--(12,0)--(11,1)--(10,0)--(11,-1)--(12,-2)--(13,-1)--(14,0)--cycle;
\draw[black] (16,-6)--(17,-5)--(18,-4)--(19,-3)--(20,-2)--(21,-1)--(20,0)--(19,1)--(18,0)--(17,1)--(16,0)--(15,1)--(14,0)--(13,-1)--(12,-2)--(13,-3)--(14,-4)--(15,-5)--(16,-6)--cycle;
\draw[black] (18,-16)--(19,-17)--(20,-18)--(21,-17)--(20,-16)--(19,-15)--(18,-14)--(19,-13)--(18,-12)--(17,-11)--(16,-10)--(15,-9)--(16,-8)--(15,-7)--(16,-6)--(15,-5)--(14,-4)--(13,-3)--(12,-2)--(11,-1)--(10,0)--(9,1)--(8,0)--(7,1)--(6,0)--(5,1)--(4,0)--(3,1)--(2,0)--(3,-1)--(4,-2)--(5,-3)--(6,-4)--(7,-5)--(8,-6)--(9,-7)--(10,-8)--(11,-9)--(12,-10)--(13,-11)--(14,-12)--(15,-13)--(16,-14)--(17,-15)--(18,-16)--cycle;
\draw[black] (21,-17)--(22,-18)--(23,-19)--(24,-18)--(25,-19)--(26,-18)--(27,-19)--(28,-18)--(29,-19)--(30,-18)--(31,-19)--(32,-18)--(33,-19)--(34,-18)--(35,-19)--(36,-18)--(37,-19)--(38,-18)--(37,-17)--(36,-16)--(35,-15)--(34,-14)--(33,-13)--(32,-12)--(31,-11)--(30,-10)--(29,-9)--(28,-8)--(27,-7)--(26,-6)--(25,-5)--(24,-4)--(23,-3)--(22,-2)--(21,-3)--(22,-4)--(21,-5)--(20,-6)--(21,-7)--(20,-8)--(19,-9)--(18,-10)--(19,-11)--(20,-12)--(19,-13)--(20,-14)--(21,-15)--(22,-16)--(21,-17)--cycle;
\draw[black] (19,-9)--(20,-8)--(21,-7)--(20,-6)--(21,-5)--(22,-4)--(21,-3)--(22,-2)--(21,-1)--(20,-2)--(19,-3)--(18,-4)--(17,-5)--(16,-6)--(17,-7)--(16,-8)--(17,-9)--(18,-10)--(19,-9)--cycle;
\draw[black] (15,-7)--(16,-8)--(17,-7)--(16,-6)--(15,-7)--cycle;
\draw[black] (17,-9)--(16,-8)--(15,-9)--(16,-10)--(17,-11)--(18,-10)--(17,-9)--cycle;
\draw[black] (17,-11)--(18,-12)--(19,-11)--(18,-10)--(17,-11)--cycle;
\draw[black] (18,-12)--(19,-13)--(20,-12)--(19,-11)--(18,-12)--cycle;
\draw[black] (18,-14)--(19,-15)--(20,-14)--(19,-13)--(18,-14)--cycle;
\draw[black] (19,-15)--(20,-16)--(21,-15)--(20,-14)--(19,-15)--cycle;
\draw[black] (20,-16)--(21,-17)--(22,-16)--(21,-15)--(20,-16)--cycle;
\node at (3,0) {$4$};
\node at (5,0) {$4$};
\node at (7,0) {$4$};
\node at (9,0) {$4$};
\node at (11,0) {$3$};
\node at (13,0) {$3$};
\node at (15,0) {$2$};
\node at (17,0) {$2$};
\node at (19,0) {$2$};
\node at (4,-1) {$4$};
\node at (6,-1) {$4$};
\node at (8,-1) {$4$};
\node at (10,-1) {$4$};
\node at (12,-1) {$3$};
\node at (14,-1) {$2$};
\node at (16,-1) {$2$};
\node at (18,-1) {$2$};
\node at (20,-1) {$2$};
\node at (5,-2) {$4$};
\node at (7,-2) {$4$};
\node at (9,-2) {$4$};
\node at (11,-2) {$4$};
\node at (13,-2) {$2$};
\node at (15,-2) {$2$};
\node at (17,-2) {$2$};
\node at (19,-2) {$2$};
\node at (21,-2) {$1$};
\node at (6,-3) {$4$};
\node at (8,-3) {$4$};
\node at (10,-3) {$4$};
\node at (12,-3) {$4$};
\node at (14,-3) {$2$};
\node at (16,-3) {$2$};
\node at (18,-3) {$2$};
\node at (20,-3) {$1$};
\node at (22,-3) {$0$};
\node at (7,-4) {$4$};
\node at (9,-4) {$4$};
\node at (11,-4) {$4$};
\node at (13,-4) {$4$};
\node at (15,-4) {$2$};
\node at (17,-4) {$2$};
\node at (19,-4) {$1$};
\node at (21,-4) {$1$};
\node at (23,-4) {$0$};
\node at (8,-5) {$4$};
\node at (10,-5) {$4$};
\node at (12,-5) {$4$};
\node at (14,-5) {$4$};
\node at (16,-5) {$2$};
\node at (18,-5) {$1$};
\node at (20,-5) {$1$};
\node at (22,-5) {$0$};
\node at (24,-5) {$0$};
\node at (9,-6) {$\mathbf 4$};
\node at (11,-6) {$\mathbf 4$};
\node at (13,-6) {$\mathbf 4$};
\node at (15,-6) {$\mathbf 4$};
\node at (17,-6) {$\mathbf 1$};
\node at (19,-6) {$\mathbf 1$};
\node at (21,-6) {$\mathbf 0$};
\node at (23,-6) {$\mathbf 0$};
\node at (25,-6) {$\mathbf 0$};
\node at (10,-7) {$4$};
\node at (12,-7) {$4$};
\node at (14,-7) {$4$};
\node at (16,-7) {$2$};
\node at (18,-7) {$1$};
\node at (20,-7) {$1$};
\node at (22,-7) {$0$};
\node at (24,-7) {$0$};
\node at (26,-7) {$0$};
\node at (11,-8) {$4$};
\node at (13,-8) {$4$};
\node at (15,-8) {$4$};
\node at (17,-8) {$1$};
\node at (19,-8) {$1$};
\node at (21,-8) {$0$};
\node at (23,-8) {$0$};
\node at (25,-8) {$0$};
\node at (27,-8) {$0$};
\node at (12,-9) {$4$};
\node at (14,-9) {$4$};
\node at (16,-9) {$3$};
\node at (18,-9) {$1$};
\node at (20,-9) {$0$};
\node at (22,-9) {$0$};
\node at (24,-9) {$0$};
\node at (26,-9) {$0$};
\node at (28,-9) {$0$};
\node at (13,-10) {$4$};
\node at (15,-10) {$4$};
\node at (17,-10) {$3$};
\node at (19,-10) {$0$};
\node at (21,-10) {$0$};
\node at (23,-10) {$0$};
\node at (25,-10) {$0$};
\node at (27,-10) {$0$};
\node at (29,-10) {$0$};
\node at (14,-11) {$4$};
\node at (16,-11) {$4$};
\node at (18,-11) {$2$};
\node at (20,-11) {$0$};
\node at (22,-11) {$0$};
\node at (24,-11) {$0$};
\node at (26,-11) {$0$};
\node at (28,-11) {$0$};
\node at (30,-11) {$0$};
\node at (15,-12) {$4$};
\node at (17,-12) {$4$};
\node at (19,-12) {$1$};
\node at (21,-12) {$0$};
\node at (23,-12) {$0$};
\node at (25,-12) {$0$};
\node at (27,-12) {$0$};
\node at (29,-12) {$0$};
\node at (31,-12) {$0$};
\node at (16,-13) {$4$};
\node at (18,-13) {$4$};
\node at (20,-13) {$0$};
\node at (22,-13) {$0$};
\node at (24,-13) {$0$};
\node at (26,-13) {$0$};
\node at (28,-13) {$0$};
\node at (30,-13) {$0$};
\node at (32,-13) {$0$};
\node at (17,-14) {$4$};
\node at (19,-14) {$3$};
\node at (21,-14) {$0$};
\node at (23,-14) {$0$};
\node at (25,-14) {$0$};
\node at (27,-14) {$0$};
\node at (29,-14) {$0$};
\node at (31,-14) {$0$};
\node at (33,-14) {$0$};
\node at (18,-15) {$4$};
\node at (20,-15) {$2$};
\node at (22,-15) {$0$};
\node at (24,-15) {$0$};
\node at (26,-15) {$0$};
\node at (28,-15) {$0$};
\node at (30,-15) {$0$};
\node at (32,-15) {$0$};
\node at (34,-15) {$0$};
\node at (19,-16) {$4$};
\node at (21,-16) {$1$};
\node at (23,-16) {$0$};
\node at (25,-16) {$0$};
\node at (27,-16) {$0$};
\node at (29,-16) {$0$};
\node at (31,-16) {$0$};
\node at (33,-16) {$0$};
\node at (35,-16) {$0$};
\node at (20,-17) {$4$};
\node at (22,-17) {$0$};
\node at (24,-17) {$0$};
\node at (26,-17) {$0$};
\node at (28,-17) {$0$};
\node at (30,-17) {$0$};
\node at (32,-17) {$0$};
\node at (34,-17) {$0$};
\node at (36,-17) {$0$};
\node at (21,-18) {$2$};
\node at (23,-18) {$0$};
\node at (25,-18) {$0$};
\node at (27,-18) {$0$};
\node at (29,-18) {$0$};
\node at (31,-18) {$0$};
\node at (33,-18) {$0$};
\node at (35,-18) {$0$};
\node at (37,-18) {$0$};
\end{tikzpicture}
\hspace{-3.5cm}
\begin{tikzpicture}[scale=0.21]
\draw[black] (20,-18)--(21,-19)--(22,-18)--(21,-17)--(20,-18)--cycle;
\draw[black] (18,0)--(17,1)--(16,0)--(15,1)--(14,0)--(15,-1)--(16,-2)--(17,-3)--(18,-4)--(19,-3)--(18,-2)--(19,-1)--(20,0)--(19,1)--(18,0)--cycle;
\draw[black] (9,-5)--(10,-4)--(9,-3)--(10,-2)--(11,-1)--(10,0)--(9,1)--(8,0)--(7,1)--(6,0)--(5,1)--(4,0)--(3,1)--(2,0)--(3,-1)--(4,-2)--(5,-3)--(6,-4)--(7,-5)--(8,-6)--(9,-5)--cycle;
\draw[black] (9,-5)--(8,-6)--(9,-7)--(10,-8)--(11,-9)--(12,-10)--(13,-11)--(14,-12)--(15,-13)--(16,-14)--(17,-15)--(18,-16)--(19,-17)--(20,-18)--(21,-17)--(20,-16)--(21,-15)--(20,-14)--(19,-13)--(18,-12)--(19,-11)--(18,-10)--(17,-9)--(18,-8)--(19,-7)--(20,-6)--(19,-5)--(18,-4)--(17,-3)--(16,-2)--(15,-1)--(14,0)--(13,1)--(12,0)--(11,1)--(10,0)--(11,-1)--(10,-2)--(9,-3)--(10,-4)--(9,-5)--cycle;
\draw[black] (22,-18)--(23,-19)--(24,-18)--(25,-19)--(26,-18)--(27,-19)--(28,-18)--(29,-19)--(30,-18)--(31,-19)--(32,-18)--(33,-19)--(34,-18)--(35,-19)--(36,-18)--(37,-19)--(38,-18)--(37,-17)--(36,-16)--(35,-15)--(34,-14)--(33,-13)--(32,-12)--(31,-11)--(30,-10)--(29,-9)--(28,-8)--(27,-7)--(26,-6)--(25,-5)--(24,-4)--(23,-3)--(22,-2)--(21,-3)--(22,-4)--(21,-5)--(20,-6)--(21,-7)--(22,-8)--(23,-9)--(22,-10)--(21,-11)--(22,-12)--(23,-13)--(22,-14)--(21,-15)--(22,-16)--(23,-17)--(22,-18)--cycle;
\draw[black] (21,-5)--(22,-4)--(21,-3)--(22,-2)--(21,-1)--(20,0)--(19,-1)--(18,-2)--(19,-3)--(18,-4)--(19,-5)--(20,-6)--(21,-5)--cycle;
\draw[black] (19,-13)--(20,-14)--(21,-15)--(22,-14)--(23,-13)--(22,-12)--(21,-11)--(22,-10)--(23,-9)--(22,-8)--(21,-7)--(20,-6)--(19,-7)--(18,-8)--(17,-9)--(18,-10)--(19,-11)--(20,-12)--(19,-13)--cycle;
\draw[black] (18,-12)--(19,-13)--(20,-12)--(19,-11)--(18,-12)--cycle;
\draw[black] (20,-16)--(21,-17)--(22,-16)--(21,-15)--(20,-16)--cycle;
\draw[black] (21,-17)--(22,-18)--(23,-17)--(22,-16)--(21,-17)--cycle;
\node at (3,0) {$4$};
\node at (5,0) {$4$};
\node at (7,0) {$4$};
\node at (9,0) {$4$};
\node at (11,0) {$3$};
\node at (13,0) {$3$};
\node at (15,0) {$2$};
\node at (17,0) {$2$};
\node at (19,0) {$2$};
\node at (4,-1) {$4$};
\node at (6,-1) {$4$};
\node at (8,-1) {$4$};
\node at (10,-1) {$4$};
\node at (12,-1) {$3$};
\node at (14,-1) {$3$};
\node at (16,-1) {$2$};
\node at (18,-1) {$2$};
\node at (20,-1) {$1$};
\node at (5,-2) {$4$};
\node at (7,-2) {$4$};
\node at (9,-2) {$4$};
\node at (11,-2) {$3$};
\node at (13,-2) {$3$};
\node at (15,-2) {$3$};
\node at (17,-2) {$2$};
\node at (19,-2) {$1$};
\node at (21,-2) {$1$};
\node at (6,-3) {$4$};
\node at (8,-3) {$4$};
\node at (10,-3) {$3$};
\node at (12,-3) {$3$};
\node at (14,-3) {$3$};
\node at (16,-3) {$3$};
\node at (18,-3) {$2$};
\node at (20,-3) {$1$};
\node at (22,-3) {$0$};
\node at (7,-4) {$4$};
\node at (9,-4) {$4$};
\node at (11,-4) {$3$};
\node at (13,-4) {$3$};
\node at (15,-4) {$3$};
\node at (17,-4) {$3$};
\node at (19,-4) {$1$};
\node at (21,-4) {$1$};
\node at (23,-4) {$0$};
\node at (8,-5) {$4$};
\node at (10,-5) {$3$};
\node at (12,-5) {$3$};
\node at (14,-5) {$3$};
\node at (16,-5) {$3$};
\node at (18,-5) {$3$};
\node at (20,-5) {$1$};
\node at (22,-5) {$0$};
\node at (24,-5) {$0$};
\node at (9,-6) {$\mathbf 3$};
\node at (11,-6) {$\mathbf 3$};
\node at (13,-6) {$\mathbf 3$};
\node at (15,-6) {$\mathbf 3$};
\node at (17,-6) {$\mathbf 3$};
\node at (19,-6) {$\mathbf 3$};
\node at (21,-6) {$\mathbf 0$};
\node at (23,-6) {$\mathbf 0$};
\node at (25,-6) {$\mathbf 0$};
\node at (10,-7) {$3$};
\node at (12,-7) {$3$};
\node at (14,-7) {$3$};
\node at (16,-7) {$3$};
\node at (18,-7) {$3$};
\node at (20,-7) {$1$};
\node at (22,-7) {$0$};
\node at (24,-7) {$0$};
\node at (26,-7) {$0$};
\node at (11,-8) {$3$};
\node at (13,-8) {$3$};
\node at (15,-8) {$3$};
\node at (17,-8) {$3$};
\node at (19,-8) {$1$};
\node at (21,-8) {$1$};
\node at (23,-8) {$0$};
\node at (25,-8) {$0$};
\node at (27,-8) {$0$};
\node at (12,-9) {$3$};
\node at (14,-9) {$3$};
\node at (16,-9) {$3$};
\node at (18,-9) {$1$};
\node at (20,-9) {$1$};
\node at (22,-9) {$1$};
\node at (24,-9) {$0$};
\node at (26,-9) {$0$};
\node at (28,-9) {$0$};
\node at (13,-10) {$3$};
\node at (15,-10) {$3$};
\node at (17,-10) {$3$};
\node at (19,-10) {$1$};
\node at (21,-10) {$1$};
\node at (23,-10) {$0$};
\node at (25,-10) {$0$};
\node at (27,-10) {$0$};
\node at (29,-10) {$0$};
\node at (14,-11) {$3$};
\node at (16,-11) {$3$};
\node at (18,-11) {$3$};
\node at (20,-11) {$1$};
\node at (22,-11) {$0$};
\node at (24,-11) {$0$};
\node at (26,-11) {$0$};
\node at (28,-11) {$0$};
\node at (30,-11) {$0$};
\node at (15,-12) {$3$};
\node at (17,-12) {$3$};
\node at (19,-12) {$2$};
\node at (21,-12) {$1$};
\node at (23,-12) {$0$};
\node at (25,-12) {$0$};
\node at (27,-12) {$0$};
\node at (29,-12) {$0$};
\node at (31,-12) {$0$};
\node at (16,-13) {$3$};
\node at (18,-13) {$3$};
\node at (20,-13) {$1$};
\node at (22,-13) {$1$};
\node at (24,-13) {$0$};
\node at (26,-13) {$0$};
\node at (28,-13) {$0$};
\node at (30,-13) {$0$};
\node at (32,-13) {$0$};
\node at (17,-14) {$3$};
\node at (19,-14) {$3$};
\node at (21,-14) {$1$};
\node at (23,-14) {$0$};
\node at (25,-14) {$0$};
\node at (27,-14) {$0$};
\node at (29,-14) {$0$};
\node at (31,-14) {$0$};
\node at (33,-14) {$0$};
\node at (18,-15) {$3$};
\node at (20,-15) {$3$};
\node at (22,-15) {$0$};
\node at (24,-15) {$0$};
\node at (26,-15) {$0$};
\node at (28,-15) {$0$};
\node at (30,-15) {$0$};
\node at (32,-15) {$0$};
\node at (34,-15) {$0$};
\node at (19,-16) {$3$};
\node at (21,-16) {$2$};
\node at (23,-16) {$0$};
\node at (25,-16) {$0$};
\node at (27,-16) {$0$};
\node at (29,-16) {$0$};
\node at (31,-16) {$0$};
\node at (33,-16) {$0$};
\node at (35,-16) {$0$};
\node at (20,-17) {$3$};
\node at (22,-17) {$1$};
\node at (24,-17) {$0$};
\node at (26,-17) {$0$};
\node at (28,-17) {$0$};
\node at (30,-17) {$0$};
\node at (32,-17) {$0$};
\node at (34,-17) {$0$};
\node at (36,-17) {$0$};
\node at (21,-18) {$2$};
\node at (23,-18) {$0$};
\node at (25,-18) {$0$};
\node at (27,-18) {$0$};
\node at (29,-18) {$0$};
\node at (31,-18) {$0$};
\node at (33,-18) {$0$};
\node at (35,-18) {$0$};
\node at (37,-18) {$0$};
\end{tikzpicture}
\]
and
\[
\begin{tikzpicture}[scale=0.21]
\draw[black] (20,-18)--(21,-19)--(22,-18)--(21,-17)--(20,-18)--cycle;
\draw[black] (8,-6)--(9,-5)--(10,-4)--(11,-3)--(12,-2)--(11,-1)--(10,0)--(9,1)--(8,0)--(7,1)--(6,0)--(5,1)--(4,0)--(3,1)--(2,0)--(3,-1)--(4,-2)--(5,-3)--(6,-4)--(7,-5)--(8,-6)--cycle;
\draw[black] (19,-9)--(20,-8)--(21,-7)--(22,-6)--(21,-5)--(20,-4)--(19,-3)--(18,-2)--(19,-1)--(20,0)--(19,1)--(18,0)--(17,1)--(16,0)--(15,1)--(14,0)--(15,-1)--(16,-2)--(17,-3)--(18,-4)--(17,-5)--(16,-6)--(17,-7)--(16,-8)--(15,-9)--(16,-10)--(17,-11)--(18,-10)--(19,-9)--cycle;
\draw[black] (17,-3)--(16,-2)--(15,-1)--(14,0)--(13,1)--(12,0)--(11,1)--(10,0)--(11,-1)--(12,-2)--(11,-3)--(10,-4)--(9,-5)--(8,-6)--(9,-7)--(10,-8)--(11,-9)--(12,-10)--(13,-11)--(14,-12)--(15,-13)--(16,-14)--(17,-15)--(18,-16)--(19,-17)--(20,-18)--(21,-17)--(20,-16)--(19,-15)--(18,-14)--(19,-13)--(18,-12)--(17,-11)--(16,-10)--(15,-9)--(16,-8)--(17,-7)--(16,-6)--(17,-5)--(18,-4)--(17,-3)--cycle;
\draw[black] (22,-18)--(23,-19)--(24,-18)--(25,-19)--(26,-18)--(27,-19)--(28,-18)--(29,-19)--(30,-18)--(31,-19)--(32,-18)--(33,-19)--(34,-18)--(35,-19)--(36,-18)--(37,-19)--(38,-18)--(37,-17)--(36,-16)--(35,-15)--(34,-14)--(33,-13)--(32,-12)--(31,-11)--(30,-10)--(29,-9)--(28,-8)--(27,-7)--(26,-6)--(25,-5)--(24,-4)--(23,-3)--(22,-2)--(21,-1)--(20,-2)--(19,-3)--(20,-4)--(21,-5)--(22,-6)--(21,-7)--(22,-8)--(23,-9)--(24,-10)--(25,-11)--(24,-12)--(23,-13)--(24,-14)--(25,-15)--(24,-16)--(23,-17)--(22,-18)--cycle;
\draw[black] (18,-2)--(19,-3)--(20,-2)--(21,-1)--(20,0)--(19,-1)--(18,-2)--cycle;
\draw[black] (17,-11)--(18,-12)--(19,-13)--(20,-14)--(19,-15)--(20,-16)--(21,-17)--(22,-18)--(23,-17)--(24,-16)--(25,-15)--(24,-14)--(23,-13)--(24,-12)--(25,-11)--(24,-10)--(23,-9)--(22,-8)--(21,-7)--(20,-8)--(19,-9)--(18,-10)--(17,-11)--cycle;
\draw[black] (18,-14)--(19,-15)--(20,-14)--(19,-13)--(18,-14)--cycle;
\node at (3,0) {$4$};
\node at (5,0) {$4$};
\node at (7,0) {$4$};
\node at (9,0) {$4$};
\node at (11,0) {$3$};
\node at (13,0) {$3$};
\node at (15,0) {$2$};
\node at (17,0) {$2$};
\node at (19,0) {$2$};
\node at (4,-1) {$4$};
\node at (6,-1) {$4$};
\node at (8,-1) {$4$};
\node at (10,-1) {$4$};
\node at (12,-1) {$3$};
\node at (14,-1) {$3$};
\node at (16,-1) {$2$};
\node at (18,-1) {$2$};
\node at (20,-1) {$1$};
\node at (5,-2) {$4$};
\node at (7,-2) {$4$};
\node at (9,-2) {$4$};
\node at (11,-2) {$4$};
\node at (13,-2) {$3$};
\node at (15,-2) {$3$};
\node at (17,-2) {$2$};
\node at (19,-2) {$1$};
\node at (21,-2) {$0$};
\node at (6,-3) {$4$};
\node at (8,-3) {$4$};
\node at (10,-3) {$4$};
\node at (12,-3) {$3$};
\node at (14,-3) {$3$};
\node at (16,-3) {$3$};
\node at (18,-3) {$2$};
\node at (20,-3) {$0$};
\node at (22,-3) {$0$};
\node at (7,-4) {$4$};
\node at (9,-4) {$4$};
\node at (11,-4) {$3$};
\node at (13,-4) {$3$};
\node at (15,-4) {$3$};
\node at (17,-4) {$3$};
\node at (19,-4) {$2$};
\node at (21,-4) {$0$};
\node at (23,-4) {$0$};
\node at (8,-5) {$4$};
\node at (10,-5) {$3$};
\node at (12,-5) {$3$};
\node at (14,-5) {$3$};
\node at (16,-5) {$3$};
\node at (18,-5) {$2$};
\node at (20,-5) {$2$};
\node at (22,-5) {$0$};
\node at (24,-5) {$0$};
\node at (9,-6) {$\mathbf 3$};
\node at (11,-6) {$\mathbf 3$};
\node at (13,-6) {$\mathbf 3$};
\node at (15,-6) {$\mathbf 3$};
\node at (17,-6) {$\mathbf 2$};
\node at (19,-6) {$\mathbf 2$};
\node at (21,-6) {$\mathbf 2$};
\node at (23,-6) {$\mathbf 0$};
\node at (25,-6) {$\mathbf 0$};
\node at (10,-7) {$3$};
\node at (12,-7) {$3$};
\node at (14,-7) {$3$};
\node at (16,-7) {$3$};
\node at (18,-7) {$2$};
\node at (20,-7) {$2$};
\node at (22,-7) {$0$};
\node at (24,-7) {$0$};
\node at (26,-7) {$0$};
\node at (11,-8) {$3$};
\node at (13,-8) {$3$};
\node at (15,-8) {$3$};
\node at (17,-8) {$2$};
\node at (19,-8) {$2$};
\node at (21,-8) {$1$};
\node at (23,-8) {$0$};
\node at (25,-8) {$0$};
\node at (27,-8) {$0$};
\node at (12,-9) {$3$};
\node at (14,-9) {$3$};
\node at (16,-9) {$2$};
\node at (18,-9) {$2$};
\node at (20,-9) {$1$};
\node at (22,-9) {$1$};
\node at (24,-9) {$0$};
\node at (26,-9) {$0$};
\node at (28,-9) {$0$};
\node at (13,-10) {$3$};
\node at (15,-10) {$3$};
\node at (17,-10) {$2$};
\node at (19,-10) {$1$};
\node at (21,-10) {$1$};
\node at (23,-10) {$1$};
\node at (25,-10) {$0$};
\node at (27,-10) {$0$};
\node at (29,-10) {$0$};
\node at (14,-11) {$3$};
\node at (16,-11) {$3$};
\node at (18,-11) {$1$};
\node at (20,-11) {$1$};
\node at (22,-11) {$1$};
\node at (24,-11) {$1$};
\node at (26,-11) {$0$};
\node at (28,-11) {$0$};
\node at (30,-11) {$0$};
\node at (15,-12) {$3$};
\node at (17,-12) {$3$};
\node at (19,-12) {$1$};
\node at (21,-12) {$1$};
\node at (23,-12) {$1$};
\node at (25,-12) {$0$};
\node at (27,-12) {$0$};
\node at (29,-12) {$0$};
\node at (31,-12) {$0$};
\node at (16,-13) {$3$};
\node at (18,-13) {$3$};
\node at (20,-13) {$1$};
\node at (22,-13) {$1$};
\node at (24,-13) {$0$};
\node at (26,-13) {$0$};
\node at (28,-13) {$0$};
\node at (30,-13) {$0$};
\node at (32,-13) {$0$};
\node at (17,-14) {$3$};
\node at (19,-14) {$2$};
\node at (21,-14) {$1$};
\node at (23,-14) {$1$};
\node at (25,-14) {$0$};
\node at (27,-14) {$0$};
\node at (29,-14) {$0$};
\node at (31,-14) {$0$};
\node at (33,-14) {$0$};
\node at (18,-15) {$3$};
\node at (20,-15) {$1$};
\node at (22,-15) {$1$};
\node at (24,-15) {$1$};
\node at (26,-15) {$0$};
\node at (28,-15) {$0$};
\node at (30,-15) {$0$};
\node at (32,-15) {$0$};
\node at (34,-15) {$0$};
\node at (19,-16) {$3$};
\node at (21,-16) {$1$};
\node at (23,-16) {$1$};
\node at (25,-16) {$0$};
\node at (27,-16) {$0$};
\node at (29,-16) {$0$};
\node at (31,-16) {$0$};
\node at (33,-16) {$0$};
\node at (35,-16) {$0$};
\node at (20,-17) {$3$};
\node at (22,-17) {$1$};
\node at (24,-17) {$0$};
\node at (26,-17) {$0$};
\node at (28,-17) {$0$};
\node at (30,-17) {$0$};
\node at (32,-17) {$0$};
\node at (34,-17) {$0$};
\node at (36,-17) {$0$};
\node at (21,-18) {$2$};
\node at (23,-18) {$0$};
\node at (25,-18) {$0$};
\node at (27,-18) {$0$};
\node at (29,-18) {$0$};
\node at (31,-18) {$0$};
\node at (33,-18) {$0$};
\node at (35,-18) {$0$};
\node at (37,-18) {$0$};
\end{tikzpicture}
\hspace{-3.5cm}
\begin{tikzpicture}[scale=0.21]
\draw[black] (11,-3)--(12,-2)--(13,-1)--(14,0)--(13,1)--(12,0)--(11,1)--(10,0)--(11,-1)--(10,-2)--(9,-3)--(10,-4)--(11,-3)--cycle;
\draw[black] (9,-5)--(10,-4)--(9,-3)--(10,-2)--(11,-1)--(10,0)--(9,1)--(8,0)--(7,1)--(6,0)--(5,1)--(4,0)--(3,1)--(2,0)--(3,-1)--(4,-2)--(5,-3)--(6,-4)--(7,-5)--(8,-6)--(9,-5)--cycle;
\draw[black] (22,-18)--(23,-19)--(24,-18)--(25,-19)--(26,-18)--(27,-19)--(28,-18)--(29,-19)--(30,-18)--(31,-19)--(32,-18)--(33,-19)--(34,-18)--(35,-19)--(36,-18)--(37,-19)--(38,-18)--(37,-17)--(36,-16)--(35,-15)--(34,-14)--(33,-13)--(32,-12)--(31,-11)--(30,-10)--(29,-9)--(28,-8)--(27,-7)--(26,-6)--(25,-7)--(24,-8)--(23,-9)--(24,-10)--(25,-11)--(24,-12)--(23,-13)--(24,-14)--(25,-15)--(24,-16)--(23,-17)--(22,-18)--cycle;
\draw[black] (18,-16)--(19,-17)--(20,-18)--(21,-19)--(22,-18)--(23,-17)--(22,-16)--(21,-15)--(22,-14)--(23,-13)--(22,-12)--(21,-11)--(22,-10)--(23,-9)--(24,-8)--(25,-7)--(26,-6)--(25,-5)--(24,-4)--(23,-3)--(22,-2)--(21,-1)--(20,0)--(19,1)--(18,0)--(17,1)--(16,0)--(15,1)--(14,0)--(13,-1)--(12,-2)--(11,-3)--(10,-4)--(9,-5)--(8,-6)--(9,-7)--(10,-8)--(11,-9)--(12,-10)--(13,-11)--(14,-12)--(15,-13)--(16,-14)--(17,-15)--(18,-16)--cycle;
\draw[black] (24,-12)--(25,-11)--(24,-10)--(23,-9)--(22,-10)--(21,-11)--(22,-12)--(23,-13)--(24,-12)--cycle;
\draw[black] (24,-16)--(25,-15)--(24,-14)--(23,-13)--(22,-14)--(21,-15)--(22,-16)--(23,-17)--(24,-16)--cycle;
\node at (3,0) {$4$};
\node at (5,0) {$4$};
\node at (7,0) {$4$};
\node at (9,0) {$4$};
\node at (11,0) {$3$};
\node at (13,0) {$3$};
\node at (15,0) {$2$};
\node at (17,0) {$2$};
\node at (19,0) {$2$};
\node at (4,-1) {$4$};
\node at (6,-1) {$4$};
\node at (8,-1) {$4$};
\node at (10,-1) {$4$};
\node at (12,-1) {$3$};
\node at (14,-1) {$2$};
\node at (16,-1) {$2$};
\node at (18,-1) {$2$};
\node at (20,-1) {$2$};
\node at (5,-2) {$4$};
\node at (7,-2) {$4$};
\node at (9,-2) {$4$};
\node at (11,-2) {$3$};
\node at (13,-2) {$2$};
\node at (15,-2) {$2$};
\node at (17,-2) {$2$};
\node at (19,-2) {$2$};
\node at (21,-2) {$2$};
\node at (6,-3) {$4$};
\node at (8,-3) {$4$};
\node at (10,-3) {$3$};
\node at (12,-3) {$2$};
\node at (14,-3) {$2$};
\node at (16,-3) {$2$};
\node at (18,-3) {$2$};
\node at (20,-3) {$2$};
\node at (22,-3) {$2$};
\node at (7,-4) {$4$};
\node at (9,-4) {$4$};
\node at (11,-4) {$2$};
\node at (13,-4) {$2$};
\node at (15,-4) {$2$};
\node at (17,-4) {$2$};
\node at (19,-4) {$2$};
\node at (21,-4) {$2$};
\node at (23,-4) {$2$};
\node at (8,-5) {$4$};
\node at (10,-5) {$2$};
\node at (12,-5) {$2$};
\node at (14,-5) {$2$};
\node at (16,-5) {$2$};
\node at (18,-5) {$2$};
\node at (20,-5) {$2$};
\node at (22,-5) {$2$};
\node at (24,-5) {$2$};
\node at (9,-6) {$\mathbf 2$};
\node at (11,-6) {$\mathbf 2$};
\node at (13,-6) {$\mathbf 2$};
\node at (15,-6) {$\mathbf 2$};
\node at (17,-6) {$\mathbf 2$};
\node at (19,-6) {$\mathbf 2$};
\node at (21,-6) {$\mathbf 2$};
\node at (23,-6) {$\mathbf 2$};
\node at (25,-6) {$\mathbf 2$};
\node at (10,-7) {$2$};
\node at (12,-7) {$2$};
\node at (14,-7) {$2$};
\node at (16,-7) {$2$};
\node at (18,-7) {$2$};
\node at (20,-7) {$2$};
\node at (22,-7) {$2$};
\node at (24,-7) {$2$};
\node at (26,-7) {$0$};
\node at (11,-8) {$2$};
\node at (13,-8) {$2$};
\node at (15,-8) {$2$};
\node at (17,-8) {$2$};
\node at (19,-8) {$2$};
\node at (21,-8) {$2$};
\node at (23,-8) {$2$};
\node at (25,-8) {$0$};
\node at (27,-8) {$0$};
\node at (12,-9) {$2$};
\node at (14,-9) {$2$};
\node at (16,-9) {$2$};
\node at (18,-9) {$2$};
\node at (20,-9) {$2$};
\node at (22,-9) {$2$};
\node at (24,-9) {$0$};
\node at (26,-9) {$0$};
\node at (28,-9) {$0$};
\node at (13,-10) {$2$};
\node at (15,-10) {$2$};
\node at (17,-10) {$2$};
\node at (19,-10) {$2$};
\node at (21,-10) {$2$};
\node at (23,-10) {$1$};
\node at (25,-10) {$0$};
\node at (27,-10) {$0$};
\node at (29,-10) {$0$};
\node at (14,-11) {$2$};
\node at (16,-11) {$2$};
\node at (18,-11) {$2$};
\node at (20,-11) {$2$};
\node at (22,-11) {$1$};
\node at (24,-11) {$1$};
\node at (26,-11) {$0$};
\node at (28,-11) {$0$};
\node at (30,-11) {$0$};
\node at (15,-12) {$2$};
\node at (17,-12) {$2$};
\node at (19,-12) {$2$};
\node at (21,-12) {$2$};
\node at (23,-12) {$1$};
\node at (25,-12) {$0$};
\node at (27,-12) {$0$};
\node at (29,-12) {$0$};
\node at (31,-12) {$0$};
\node at (16,-13) {$2$};
\node at (18,-13) {$2$};
\node at (20,-13) {$2$};
\node at (22,-13) {$2$};
\node at (24,-13) {$0$};
\node at (26,-13) {$0$};
\node at (28,-13) {$0$};
\node at (30,-13) {$0$};
\node at (32,-13) {$0$};
\node at (17,-14) {$2$};
\node at (19,-14) {$2$};
\node at (21,-14) {$2$};
\node at (23,-14) {$1$};
\node at (25,-14) {$0$};
\node at (27,-14) {$0$};
\node at (29,-14) {$0$};
\node at (31,-14) {$0$};
\node at (33,-14) {$0$};
\node at (18,-15) {$2$};
\node at (20,-15) {$2$};
\node at (22,-15) {$1$};
\node at (24,-15) {$1$};
\node at (26,-15) {$0$};
\node at (28,-15) {$0$};
\node at (30,-15) {$0$};
\node at (32,-15) {$0$};
\node at (34,-15) {$0$};
\node at (19,-16) {$2$};
\node at (21,-16) {$2$};
\node at (23,-16) {$1$};
\node at (25,-16) {$0$};
\node at (27,-16) {$0$};
\node at (29,-16) {$0$};
\node at (31,-16) {$0$};
\node at (33,-16) {$0$};
\node at (35,-16) {$0$};
\node at (20,-17) {$2$};
\node at (22,-17) {$2$};
\node at (24,-17) {$0$};
\node at (26,-17) {$0$};
\node at (28,-17) {$0$};
\node at (30,-17) {$0$};
\node at (32,-17) {$0$};
\node at (34,-17) {$0$};
\node at (36,-17) {$0$};
\node at (21,-18) {$2$};
\node at (23,-18) {$0$};
\node at (25,-18) {$0$};
\node at (27,-18) {$0$};
\node at (29,-18) {$0$};
\node at (31,-18) {$0$};
\node at (33,-18) {$0$};
\node at (35,-18) {$0$};
\node at (37,-18) {$0$};
\end{tikzpicture}
\]
We cannot write $G$ as as sum of two integral GT-patterns in $\polyG^\ast_{\lambda/\mu,w}$,
since that would, in particular, express the partition in row $15$ (entries in bold) as a sum of two only 
two integer partitions with size $18$ and $9$ parts. 
This is impossible due to \cref{ex:partitionpolytope}.
\end{example}
Here is how the full example was constructed. 
The starting point is \cref{ex:partitionpolytope},
which guarantees that $G$ is not a sum of two GT-patterns in $\polyG^\ast_{\lambda/\mu,w}$.
The idea is to extend the four partitions in \cref{eq:4partitions} to GT-patterns. 
Note that there are several restrictions:
the inequalities in \cref{eq:gtinequalities} have to be satisfied and all four GT-patterns must
have the same parameters $\lambda/\mu$ and $w$.
Finally, $G_1 + \dots + G_4$ has to result in a GT-pattern with only even entries.

With these restrictions, starting from \cref{ex:partitionpolytope} 
and extending the partitions row by row by a computer search,
eventually leads to the above example.
Note that this approach is not guaranteed to work in general ---
using the other point, $(6,6,5,5,5,4,2,2,1)$ in $2\cdot \polyP_{18,9}$ with a similar setup 
is impossible to extend to a full example.

\medskip 
\begin{remark}
It is known (see \cite{Rassart2004}) that the polytope $\polyG_{\lambda/\mu,w}$ 
always has a polynomial Ehrhart functions, despite not being 
integral in general. Lattice points in $\polyG_{\lambda/\mu,w}$ are enumerated by the skew \emph{Kostka coefficient},
$K_{\lambda/\mu,w}$. It is conjectured in \cite{King04stretched} that the 
function $n \to K_{n\lambda/n\mu,nw}$ is a polynomial with non-negative coefficients.
Equivalently, the conjecture states that $\polyG_{\lambda/\mu,w}$ has an Ehrhart polynomial with non-negative coefficients.

In the light of \cref{ex:gtface}, it might be that eventual counterexamples
to the Kostka coefficient conjecture only appear in extremely ($\geq 200$) high dimensions.
\end{remark}

\section{Hooks, volume and lattice point enumeration}

We now leave the Gelfand--Tsetlin polytopes and focus on certain order polytopes.

\medskip

Let $P$ be a partial order on $[n]$. A \defin{linear extension} of $P$ is a permutation $\pi$
such that $\pi(i) < \pi(j)$ if $i <_P j$. Counting linear extensions of posets 
is $\#P$-hard in general, \cite{Brightwell1991} but for some posets there are efficient formulas.

Given a partition $\lambda$, we construct a poset $P_\lambda$ presenting its Hasse diagram 
as the Ferrers diagram of shape $\lambda$ rotated $135^\circ$,
see \cref{ex:youngdiagrams} below where covering relations for two such posets are illustrated.
Linear extensions of $P_\lambda$ correspond to standard Young tableaux of shape $\lambda$.
Since $P_\lambda$ defines an order polytope, it follows from \cite{Stanley86TwoPosetPolytopes}
that the number of linear extensions of $P_\lambda$ is the same as the normalized volume of $\oPoly(P_\lambda)$.

We know that the number of linear extensions of $P_\lambda$ can be computed using 
the \emph{Hook formula}, \cite{FrameRobinsonThrall}: 
\begin{theorem}[Hook formula]
\[
\#\{\text{linear extensions of $P_\lambda$}\} = \frac{n!}{\prod_{ (i,j) \in \lambda} \hook_\lambda((i,j))}.
\]
\end{theorem}
We define the hook $\hook_\lambda( (i,j) )$ as $ \lambda_i + \lambda'_j - i - j +1 $.
The hooks appear as vertex labels in the diagrams in \cref{ex:youngdiagrams}.

Furthermore, the volume of a polytope is determined by the leading coefficient in its Ehrhart polynomial.
It is therefore natural to ask if we can determine the entire Ehrhart polynomial of $\oPoly( P_\lambda )$
from the multiset of hooks of the diagram $\lambda$.

\begin{example}\label{ex:youngdiagrams}
Consider the diagrams $\lambda=(8,5,4)$ and $\mu=(7,7,2,1)$.
These diagrams have the same multiset of hooks (indicated by vertex labels) and thus the volume
of $\oPoly( P_\lambda )$ and $\oPoly( P_\mu )$ coincide. 
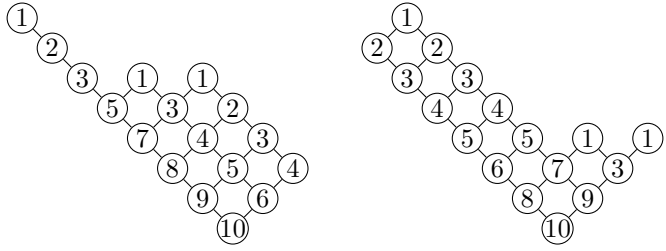
\[
\begin{tikzpicture}[xscale=0.4,yscale=0.4]
\tikzset{main node/.style={circle,draw,minimum size=0.4cm,inner sep=0pt}}
    \node[main node] (1) at ( 0, 0) {10};
    \node[main node] (2) at ( 1, 1) {6};
    \node[main node] (3) at ( 2, 2) {4};

    \node[main node] (4) at ( -1, 1) {9};
    \node[main node] (5) at ( 0, 2) {5};
    \node[main node] (6) at ( 1, 3) {3};

    \node[main node] (7) at ( -2, 2) {8};
    \node[main node] (8) at ( -1, 3) {4};
    \node[main node] (9) at (  0, 4) {2};

    \node[main node] (10) at ( -3, 3) {7};
    \node[main node] (11) at ( -2, 4) {3};
    \node[main node] (12) at ( -1, 5) {1};

    \node[main node] (13) at ( -4, 4) {5};
    \node[main node] (14) at ( -3, 5) {1};

    \node[main node] (15) at ( -5, 5) {3};
    \node[main node] (16) at ( -6, 6) {2};
    \node[main node] (17) at ( -7, 7) {1};
 
    \draw[black] (1) -- (2) -- (3);
    \draw[black] (4) -- (5) -- (6);
    \draw[black] (7) -- (8) -- (9);
    \draw[black] (10) -- (11) -- (12);
    \draw[black] (13) -- (14);
    \draw[black] (1)--(4)--(7)--(10)--(13)--(15)--(16)--(17);
    \draw[black] (2)--(5)--(8)--(11)--(14);
    \draw[black] (3)--(6)--(9)--(12);
\end{tikzpicture}
\qquad
\begin{tikzpicture}[xscale=0.4,yscale=0.4]
\tikzset{main node/.style={circle,draw,minimum size=0.4cm,inner sep=0pt}}
    \node[main node] (1) at ( 0, 0) {10};
    \node[main node] (2) at ( 1, 1) {9};
    \node[main node] (3) at ( 2, 2) {3};
    \node[main node] (4) at ( 3, 3) {1};

    \node[main node] (5) at ( -1, 1) {8};
    \node[main node] (6) at ( 0, 2) {7};
    \node[main node] (7) at ( 1, 3) {1};

    \node[main node] (8) at ( -2, 2) {6};
    \node[main node] (9) at ( -1, 3) {5};

    \node[main node] (10) at ( -3, 3) {5};
    \node[main node] (11) at ( -2, 4) {4};

    \node[main node] (12) at ( -4, 4) {4};
    \node[main node] (13) at ( -3, 5) {3};

    \node[main node] (14) at ( -5, 5) {3};
    \node[main node] (15) at ( -4, 6) {2};

    \node[main node] (16) at ( -6, 6) {2};
    \node[main node] (17) at ( -5, 7) {1};

    \draw[black] (1) -- (2) -- (3) -- (4);
    \draw[black] (5) -- (6) -- (7);
    \draw[black] (8) -- (9);
    \draw[black] (10) -- (11);
    \draw[black] (12) -- (13);
    \draw[black] (14) -- (15);
    \draw[black] (16) -- (17);
    \draw[black] (1)--(5)--(8)--(10)--(12)--(14)--(16);
    \draw[black] (2)--(6)--(9)--(11)--(13)--(15)--(17);
    \draw[black] (3)--(7);
\end{tikzpicture}
\]
However, these polytopes have different Ehrhart polynomials, 
one can easily verify that $\ehr_{\oPoly( P_\lambda )}(1) = 115$
but $\ehr_{\oPoly( P_\mu )}(1) = 134$.
Using the software \texttt{Normaliz}, one can find the complete Ehrhart polynomials
\[
\ehr_{\oPoly( P_\lambda )}(k) = 7 F(k)  (k+3) (k+4) \left(4 k^3+35 k^2+101 k+90\right)
\]
\[
\ehr_{\oPoly( P_\mu )}(k) = F(k) (k+6) (k+7) \left(28 k^3+161 k^2+301 k+180\right)
\]
where
\[
 F(k) = \frac{(k+1) (k+2)^2 (k+3)^2 (k+4)^2 (k+5)^2 (k+6) (k+7) (k+8)}{36578304000}.
\]

\end{example}
One can  ask if the multiset of hooks perhaps determine other properties of $\oPoly( P_\lambda )$.
In light off \cref{ex:youngdiagrams}, it is rather surprising that for any $\lambda$, the Ehrhart polynomial of the restriction
\begin{align}\label{eq:hooksonahyperplane}
 P_\lambda \cap \{ (x_1,\dotsc,x_n) \in \setR^n :  x_1 + \dotsc + x_n = 1 \}
\end{align}
\emph{is} determined by the multiset of hooks.
This is proved in \cite[Lemma 1]{Pak2001}, which restates the \emph{Hook Content Formula}
as an equality between lattice point cardinalities.

\bigskip 

There is a similar hook formula for counting the number 
of linear extensions of trees, see \cite{Knuth1998ArtOfProgramming}.
Here, hook values are simply given by subtree size.

\begin{example}
The trees $T$ and $T'$ have the same multiset of hooks but different Ehrhart polynomials.
Each node $v$ is decorated by its hook value --- the size of the subtree with root $v$.

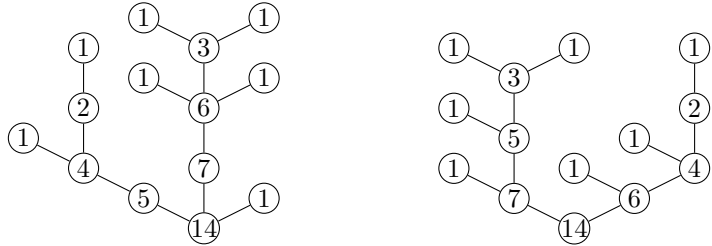
\[
 \begin{tikzpicture}[xscale=0.4,yscale=0.4]
\tikzset{main node/.style={circle,draw,minimum size=0.4cm,inner sep=0pt}}
    \node[main node] (1) at ( 0, 0) {14};
    \node[main node] (2) at ( 0, 2) {7};
    \node[main node] (3) at ( 0, 4) {6};
    \node[main node] (4) at ( 0, 6) {3};

    \node[main node] (5) at ( -2, 7) {1};
    \node[main node] (6) at ( 2, 7) {1};
    \node[main node] (7) at ( -2, 5) {1};

    \node[main node] (8) at ( 2, 5) {1};
    \node[main node] (9) at ( 2, 1) {1};

    \node[main node] (10) at ( -2, 1) {5};
    \node[main node] (11) at ( -4, 2) {4};

    \node[main node] (12) at ( -4, 4) {2};
    \node[main node] (13) at ( -4, 6) {1};

    \node[main node] (14) at ( -6, 3) {1};

    \draw[black] (1) -- (2) -- (3) -- (4) -- (5);
    \draw[black] (4) -- (6);
    \draw[black] (3) -- (7);
    \draw[black] (3) -- (8);
    \draw[black] (1) -- (9);
    \draw[black] (1) -- (10) -- (11) -- (12) -- (13);
    \draw[black] (11) -- (14);
\end{tikzpicture}
\qquad\qquad\qquad
\begin{tikzpicture}[xscale=0.4,yscale=0.4]
\tikzset{main node/.style={circle,draw,minimum size=0.4cm,inner sep=0pt}}
    \node[main node] (1) at ( 2, 1) {14};
    \node[main node] (2) at ( 0, 2) {7};
    \node[main node] (3) at ( 0, 4) {5};
    \node[main node] (4) at ( 0, 6) {3};

    \node[main node] (5) at ( -2, 7) {1};
    \node[main node] (6) at ( 2, 7) {1};
    \node[main node] (7) at ( -2, 5) {1};

    \node[main node] (8) at ( -2, 3) {1};
    \node[main node] (9) at ( 4, 2) {6};

    \node[main node] (10) at ( 2, 3) {1};
    \node[main node] (11) at ( 6, 3) {4};

    \node[main node] (12) at ( 6, 5) {2};
    \node[main node] (13) at ( 6, 7) {1};

    \node[main node] (14) at ( 4, 4) {1};

    \draw[black] (1) -- (2) -- (3) -- (4) -- (5);
    \draw[black] (4) -- (6);
    \draw[black] (3) -- (7);
    \draw[black] (2) -- (8);
    \draw[black] (1) -- (9) -- (10);
    \draw[black] (9) -- (11) -- (12) -- (13);
    \draw[black] (11) -- (14);
\end{tikzpicture}
\]

We have that
\begin{align*}
\ehr_{\oPoly(T)}(n) &= F(n) (
51480 n^9+1182984 n^8+11490017 n^7 +61564083 n^6 \\
  &+199510913 n^5 +404186041 n^4+512043278 n^3 \\
  &+393196652 n^2+167403432 n+30270240)
\end{align*}
and
\begin{align*}
\ehr_{\oPoly(T')}(n) &= F(n)(51480 n^9+1173975 n^8+11327855n^7 + 60383085 n^6 \\
   &+195027707 n^5 + 394660980 n^4+500753090 n^3 \\
   &+386259540n^2+165675888 n+30270240)
\end{align*}
where $F(n)=\frac{(n+1) (n+2) (n+3) (n+4) (n+5)}{3632428800}$.

Note that it is enough to notice that $\ehr_{\oPoly(T)}(1) \neq \ehr_{\oPoly(T')}(1)$ to deduce that the Ehrhart polynomials are different.
This computation can be done by brute force. To compute the full Ehrhart polynomial, more sophisticated software
such as \texttt{lattE} is required.
\end{example}

By combining the result in \cite{BjornerWachs89} and the theory of $P$-partitions by Stanley \cite{StanleyEC1},
one obtains the analogous statement of \eqref{eq:hooksonahyperplane} for trees.


\subsection*{Acknowledgement}

The author would like to thank Greta Panova for pointing out 
reference \cite{Pak2001} as well as for several insightful comments.
We also thank Winfied Bruns for the computations in \texttt{Normaliz}.

This work has been funded by the \emph{Knut and Alice Wallenberg Foundation} (2013.03.07).

\bibliographystyle{amsalpha}
\bibliography{bibliography}

\end{document}